\documentclass[12pt]{article}
\usepackage{graphicx}
\usepackage[usenames]{color}
\usepackage{natbib}

\usepackage{amssymb}
\usepackage{epsfig}
\usepackage{epsf}
\usepackage{float}
\usepackage{a4wide}

\newcommand{\mc}[1]{{\mathcal #1}}

\newcommand{\bb}[1]{{\mathbb #1}}
\newcommand{\reff}[1]{(\ref{#1})}

\newcommand{\Cal}[1]{{\mathcal #1}}

\newcommand{\argmin}{\mathop{\mathrm{arg\,min}}}

\def\emptysquare{{\hbox{\vrule height6pt width0.6pt depth0pt%
\vbox{\hrule height0.6pt width4.8pt depth0pt%
\vglue4.8pt%
\hrule height0.6pt width4.8pt depth0pt}%
\vrule height6pt width0.6pt depth0pt}}}

\def\qed{\unskip\nobreak
\hfil\penalty50\hskip1.75em\null\nobreak\hfil\emptysquare
{\parfillskip=0pt \finalhyphendemerits=0 \par}\medskip}

\newtheorem{thm}{Theorem}[section]
\newtheorem{lemma}[thm]{Lemma}

\newtheorem{definition}[thm]{Definition}

\makeatletter \@addtoreset{equation}{section} \makeatother

\renewcommand\thetable{\thesection.\@arabic\c@table}

\def\reff#1{(\ref{#1})}
\def\emph#1{{\it#1}}

\def\E{{\mathbb E}}
\def\P{{\mathbb P}}
\def\R{{\mathbb R}}

\def\eps{\varepsilon}

\def\gen{{\rm gen}}

\def\1{{1\kern-.25em\hbox{\rm I}}}
\def\one{\mathbf 1\kern-.15em}
\def\eu{{1\kern-.25em\hbox{\sm I}}}

\def\R{{\mathbb R}}  
\def\P{{\mathbb P}}  
\def\E{{\mathbb E}}  

\let\cal=\mathcal

\def\BB{{\cal B}}

\def\FF{{\cal F}}

\def\KK{{\cal K}}

\def\TT{{\cal T}}

\def\text#1{\quad{\hbox{#1}}\quad}

\def\proof{{\vskip 2mm \noindent\bf Proof\hskip 3mm }}

\def\endproof{$\square$ \vskip 2mm}
\def\square{\ifmmode\sqr\else{$\sqr$}\fi}
\def\sqr{\vcenter{
         \hrule height.1mm
         \hbox{\vrule width.1mm height2.2mm\kern2.18mm\vrule width.1mm}
         \hrule height.1mm}}                  

\begin{document}
\begin{title}
{Limit theorems for sequences of random trees}

\author{David Balding, Pablo A. Ferrari, Ricardo Fraiman\footnote{Corresponding author. Postal address:
Departamento de Matem\'atica y Ciencias, Universidad de San Andr\'es, Vito Dumas 284, 1644, Victoria, Argentina.
\it Tel\rm.: 54-11-47257062. \it Fax\rm: 54-11-47257010.  Email
address: \tt rfraiman@udesa.edu.ar\rm} \\
and Mariela Sued
\bigskip
\\ Department of Epidemiology and Public Health, Imperial
College, England
\\ Instituto de Matem\'atica e Estat\'{\i}stica, Univ. de S\~ao Paulo, Brasil
\\ Departamento de Matem\'{a}tica y Ciencias, Univ. de
San Andr\'es, Argentina \\ and Centro de Matem\'atica,
 Univ. de la Rep\'ublica, Uruguay
\\ Instituto del C\'alculo, Univ. de Buenos Aires, Argentina}

\end{title}
\maketitle

Running Title: Limit theorems on random trees

\




\begin{abstract}
  We consider a random tree and introduce a metric in the space of
  trees to define the ``mean tree'' as the tree minimizing the
  average distance to the random tree. When the resulting metric space
  is compact we have laws of large numbers and central limit theorems
  for sequence of independent identically distributed random trees.
  As application we propose tests to check if two samples of random
  trees have the same law.
\end{abstract}
\section{Introduction}
\label{s1}

Random trees have long been an important modelling tool. In
particular, trees are useful when a collection of observed objects
are all descended from a common ancestral object via a process of
duplication followed by gradual differentiation.  This
characterizes the process of natural evolution, and also any form
of information that over time is successively replicated, and
transmitted with occasional error.  There are two broad approaches
to constructing random evolutionary trees: forwards in time
``branching process'' models, such as the Galton-Watson process,
and backwards-in-time ``coalescent'' models such as Kingman's
coalescent (Kingman, 1982).

We prove law of large numbers and an invariance principle for random trees
defined in a metric space and propose a Kolmogorov-Smirnov-type goodness-of-fit
test.

Our trees have a special vertex called root and evolve forward in
time in discrete generations; each parent node (or vertex) has up
to $m$ offspring nodes in the next generation. The set of possible
vertices is called $\widetilde V$. A tree is a function $x:
\widetilde V\to\{0,1\}$, where $x(v)$ indicates if the vertex
$v\in\widetilde V$ is present in $x$, with the restriction that a
vertex cannot be present if its mother is not. Call $\TT$ the
resulting space of trees; when $\widetilde V$ is finite, $\TT$ is
also finite. In the general case, $\TT$ is a closed subset of the
compact product space $\{0,1\}^{\widetilde V}$. Since the product
topology is the one where convergence is in each coordinate, the
topology may be induced by different distances. In this setting
$\BB$, the Borel $\sigma$-field, is the same as the one generated
by the projections.  A similar setup was proposed by Otter (1949)
and Neveu (1986), see Kurata and Minami (2004). For a probability
measure $\nu$ on $\TT$, a random tree with law $\nu$ and a
distance $d$ on $\TT$, the $d$-mean related to $\nu$ is defined as
the tree (or set of trees) that minimizes the $\nu$-average
$d$-distance to the random tree. Other tree spaces and metrics are
briefly discussed in Section \ref{s8}.

We consider a sample of independent and identically distributed random elements
of a compact metric space with law $\nu$ and a unique $d$-mean. We prove that
the empiric $d$-mean of the sample converges to the $d$-mean related to $\nu$ as
the size of the sample goes to infinity. Hence the empiric $d$-mean is a
consistent estimator for the $d$-mean related to $\nu$. The result applies to
metric spaces of trees that may have infinitely many vertices.  The law of large
numbers on metric spaces with negative curvature has been addressed by Herer
(1992), de Fitte (1997) and Es-Sahib and Heinich (1999). Our space is not of
negative curvature, as shown in Section \ref{s8}. For compact metric spaces, a
strong law of large numbers have been obtained in Sverdrup-Thygenson (1981),
which is used in our setting.

We show an invariance principle for the random processes $(g_n(y)-g(y),\,y\in
\TT)$, where $g_n(y)$ is the average of the distances from $y$ to the points of
a sample of size $n$ and $g(y)$ is the average of the distances from $y$ to the
random tree with law $\nu$ from where the sample is obtained. The proof is based
on a theorem by Ledoux and Talagrand (1991); we build up a probability measure
on the space of trees that satisfies the ``majorizing measure condition'' for a
particular family of distances.

The invariance principle implies the approximate distribution of
\begin{equation}
  \label{r1}
  \max_{y\in \TT} |g_n(y)-g(y)|,
\end{equation}
is known. We propose \reff{r1} as statistic for a universal Kolmogorov-type
goodness of fit test and the analogous for the two-sample problem. In general
$(g_n(y)-g(y),\,y\in \TT)$ does not identify the measure $\nu$. Busch et al
(2006) show that $(g_n(y)-g(y),\,y\in \TT)$ identifies the vertex-marginals
$(\nu\{x:\,x(v)=1\},\,v\in\widetilde V)$ and viceversa. The vertex-marginals do
not always identify the measure but they do if the tree is constructed in a
Markovian way; examples include Galton-Watson and other related processes.

As far as we know the Otter-Neveu set-up has not been used before to construct
statistical tools for random trees. With this structure the law of large numbers
and invariance principles are quite straightforward and the statistic \reff{r1}
arises naturally to perform goodness-to-fit tests. The computation of the
statistic \reff{r1} requires in principle an exponential number of steps in the
number of possible vertices.  Busch et al (2006) show that the search of the
maximum in \reff{r1} is equivalent to the search of the minimal cut in an
associated network; a technique coming from image reconstruction. This makes the
test viable for reasonable big trees.

The critical values related to the statistic \reff{r1} depend on the
distribution $\nu$. To compute them it is usually necessary to simulate trees
with the tested distribution or to perform bootstrap. Our test has been applied
to samples of Galton-Watson related processes obtained by simulation and to a
classification of FGF protein families (Busch et al 2006). In both cases the
test has been successfull to distinguish different laws, even when the mean tree
is the same for the two samples.


In Section~\ref{s2} we introduce the space of trees as a metric space and define
the $d$-mean tree. In Section~\ref{s3} we prove the law of large numbers.  In
Section~\ref{s5} we give some examples and in Section~\ref{s6} we prove the
invariance principle.  In Section~\ref{s7} we describe the statistical
applications.  In Section \ref{s8} we show that our space is not of negative
curvature and discuss some other possible metrics.

\section{A metric space of rooted trees}
\label{s2}

Let $\widetilde V= \{1,11,12,\dots,1m,\dots\}$ the set of finite sequences of
numbers in $A=\{1,\dots,m\}$ starting with 1, with $m$ a natural number.
Elements of $\widetilde V$ are called \emph{vertices}; the vertex $1$ is called
\emph{root}. The \emph{full tree} is the oriented graph $\tilde x=(\widetilde
V,\widetilde E)$ with edges $\widetilde E\subset \widetilde V\times \widetilde
V$ given by $\widetilde E=\{(v,va)\,:\, v\in\widetilde V,\,a\in A\}$, where $va$
is the sequence obtained by juxtaposition of $v$ and $a$. In the full tree each
node or vertex has exactly $m$ outgoing edges to her offsprings and one ingoing
edge from her mother, except for the root that has no ingoing edges.  The node
$v=a_1\dots a_k$ is said to belong to the \emph{generation} $k$; in this case we
write $\gen(v)=k$. Generation 1 has only one node: the root of the tree.


We define a tree as a function $x:\widetilde V\to\{0,1\}$ satisfying,
for all $v\in \widetilde V$ and $a\in A$,
\begin{eqnarray}\label{i1}
x(v)&\ge&x(va).
\end{eqnarray}
Abusing notation, we identify $x$ with the graph $x=(V_x,E_x)$ where
\begin{eqnarray}
  \label{p1}
&&V_x=\{v\in \widetilde V\,:\, x(v)=1\},\\
&&E_x = \{(v,va) \in \widetilde E\,:\, x(v)= x(va)=1\}\;.
\end{eqnarray}
Let $\TT$ be the set of trees of this form.  Condition \reff{i1} in effect
requires that for $x\in\TT$, every node in $x$ must have a parent node in each
previous generation back to the root.

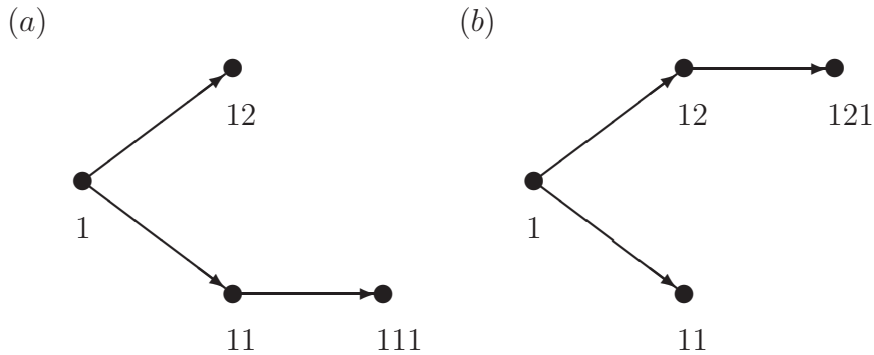
\begin{figure}[ht]
\setlength{\unitlength}{0.5mm}
\phantom{1}\hskip 2cm
\begin{picture}(200,100)
\thicklines
\put(-20,90){($a$)}
\put(0,50){\circle*{5}}\put(-2,35){1}
\put(0,50){\vector(4,3){38}}
\put(0,50){\vector(4,-3){38}}
\put(40,80){\circle*{5}}\put(38,65){12}
\put(40,20){\circle*{5}}\put(38,5){11}
\put(40,20){\vector(1,0){38}}
\put(80,20){\circle*{5}}\put(78,5){111}
\put(100,90){($b$)}
\put(120,50){\circle*{5}}\put(118,35){1}
\put(120,50){\vector(4,3){38}}
\put(120,50){\vector(4,-3){38}}
\put(160,80){\circle*{5}}\put(158,65){12}
\put(160,20){\circle*{5}}\put(158,5){11}
\put(160,80){\vector(1,0){38}}
\put(200,80){\circle*{5}}\put(198,65){121}
\end{picture}
\caption{Two finite trees both with 3 generations and 2
terminal nodes.}\label{fig1}
\end{figure}

A finite tree is characterized by the set of its terminal nodes.
For example, the trees in Figure \ref{fig1} are (a) \{111,12\},
and (b) \{11,121\}.

The product topology on $\{0,1\}^{\widetilde V}$ is the smaller
for which the projections are continuous. By projection we mean
the family of functions $\pi_v:\{0,1\}^{\widetilde V}\to \{0,1\}$
that map $x\to x(v)$. In this topology $x_n$ converges to $x$ if
and only if $x_n(v)$ converges to $x(v)$ for all $v\in \widetilde
V$. Since at each vertex we have values in $\{0,1\}$, convergence
means that for each $v$ there exists $n(v)$ such that if $n\geq
n(v)$ then $x_n(v)=x(v)$. This condition guarantees that $\TT$,
the space of trees, is a closed set in $\{0,1\}^{\widetilde V}$
and hence also compact.

As it is done in interacting particle systems (see Liggett (1985))
we consider the sigma algebra $\BB$ generated by the cylinders
$\{x\in\TT\,:\, x(v)=1\}$, $v\in\widetilde V$; this is just the
Borel sigma field generated by the product topology.

We provide $\TT$ with a \emph{distance} $d$, so that $(\TT,d)$ is a
\emph{metric space}. We use the family of distances in $\TT$ defined by
\begin{equation}
   \label{p2}
   d(x,y) = \sum_{v\in \widetilde V} |x(v)-y(v)| \phi(v)\;,
\end{equation}
for some strictly positive function $\phi:\widetilde V\to\R^+$
satisfying $\sum_{v\in \widetilde V} \phi(v)<\infty$. In this case,
the distance between the two trees of Figure \ref{fig1} is
$d(a,b)=\phi(111)+\phi(121)$.

This distance is compatible with the product
topology and hence, the notion of convergence under any of these
metrics is the same as the induced by the product topology.

Otter (1949) and Neveu (1986) propose a similar construction, but to deal with
unbounded number of offsprings, they ask each vertex $v$ and natural $a$ to
satisfy $x(v(a+1))\le x(va)$; informally, the presence of a brother in the tree
implies that all older brothers are also present. Their distance, also
compatible with the product topology, is defined by
\begin{equation}
  \label{on1}
  d_{ON}(x,t)=\exp\big(-\max \{k\,:\,x(v) = t(v)\, \hbox{ for all } v \hbox{ such
    that } \gen(v)\le k\}\big).
\end{equation}
See Kurata and
Minami (2004) for a review of those papers.

\paragraph*{\bf Random trees}
A \emph{random tree} with distribution $\nu $ is a measurable
function
\begin{equation}
  \label{p3}
  T:\Omega\to\TT\;\hbox{such that}\; \P(T\in A) = \int_A \nu(dx) \;.
\end{equation}
for any Borel set  $A\in \Cal B$, where $(\Omega,\mathcal F,\P)$ is a
probability space and $\nu$ a probability on $(\TT,\Cal B)$.

The expected distance from a tree $y$ to a random tree $T$ is defined by
\begin{eqnarray}
  \label{px10}
g(y)&:=&  \E(d(T,y))
\,=\,\int_\TT d(x,y) \;\nu(dx)\\
&=& \sum_{x\in \TT}\nu(x)d(x,y)\quad\hbox{ (in the discrete case).}
\end{eqnarray}

\begin{definition}
  The expected value or $d$-mean of a random tree $T$ is the set (of
  trees) $\E_d T$ that minimizes the expected distance to $T$:
\begin{eqnarray}
  \label{p4}
  \E_d T &:=&
\argmin_{y\in\TT}\, g(y)\,.
\end{eqnarray}
\end{definition}

The set $\E_d T$ might be empty, but if $\TT$ is compact, then $\E_d
T$ is not empty (see Section ~\ref{s3}).  Any element of the set
$\E_dT $ is also called a $d$-mean.  Since $\E_dT$ depends only on the
distribution $\nu$ induced by $T$ on $\TT$, it may also be denoted as
$\E_d(\nu)$.  The elements of $\E_d(\nu)$ are also called $d$-centers.
The notion of expected value depends on the distance $d$; in
particular, for random variables in $\R^k$ we may obtain the usual
mean, the median and the mode as illustrated in Section~\ref{s5}.

\noindent{\it Example:} In the Galton-Watson branching process the numbers of
offspring of distinct nodes are i.i.d.  In the special case that they have the
Binomial(2,$p$) distribution with $p\in [0,1]$, the offspring number is 0, 1, or
2 with probabilities $(1{-}p)^2$, $2p(1{-}p)$, and $p^2$. Letting $k_0 =
\max\{k\in \{0,1,\dots\} : p^k \ge 1/2\}$ there are two cases: (a) if $p^{k_0} >
1/2$ there is only one mean tree $x$ satisfying $x(v)=1$ if and only if
$\gen(v)< k_0$ and (b) if $p^{k_0} = 1/2$ the mean tree is the set of trees with
$x(v)=1$ if $\gen(v)< k_0$, $x(v)\in \{0,1\}$ if $\gen(v)= k_0$ and $x(v)=0$ if
$\gen(v)> k_0$.  In particular, if $p<1/2$ the mean tree is the empty tree.

Let $T_1,\dots,T_n$ be a random sample of $T$ (independent random
trees with the same law as $T$). The empiric measure associated to the
sample is denoted by $\mu_n$ and it is given by
\begin{equation}
\label{emptree}
\mu_n\;:=\;\frac1n\sum_{i=1}^n\delta
_{T_i}\;,\quad\mu_n(x)\;:=\;\frac1n\sum_{i=1}^n\one
_{\{T_i=x\}}\;,
\end{equation}
where $\delta_x$ is the point mass at $x$ and $\one_A$ is the
indicator function of the set $A$.  Associated to this measures we
define the empiric expected distance
of a tree $y$ to the sample by
\begin{equation}
  \label{px1}
  g_n(y) := \int_\TT d(x,y)\; \mu_n(dx) = \frac1n\sum_{i=1}^n
d(T_i,y)\,,
\end{equation}
and as in \reff{p4} the \emph{empiric mean} tree (empiric $d$-center,
sample $d$-mean) as the random set given by
\begin{equation}
  \label{p5}
  \overline T_n := \argmin_{y\in\TT} g_n(y)\,. 
\end{equation}
The empirical mean is like a consensus tree. If $n$ is odd, the
empirical mean is unique; it just includes all vertices that are in
more than half of the trees. If $n$ is even, it is not unique but
there is a ``shortest'' and ``largest'' empirical mean tree, and every
subtree of the largest empirical mean tree which contains the shortest
empirical mean tree is on the set of empirical $d$-means. This is a
nice property from the robustness point of view.

\section{Law of large numbers}
\label{s3}

If $\nu$ is defined on a finite set of trees the following law of large numbers
follows immediately.
\begin{thm}
  \label{p6}
  Let $(\TT_0,d)$ be a finite tree space with metric $d$. Let
  $T\in\TT_0$ be a random tree with law $\nu$ such that $\E_d T$ has only one
element (also denoted by $\E_dT$).
   Let $\{T_n,\, n\ge 1\}$ be
  an i.i.d. sequence of random trees with law $\nu$. If $y_n$ is any of the
  empiric mean trees of $\{T_1,\dots,T_n\}$, that is $y_n\in \overline T_n$, then
\begin{equation}
  \label{p7}
  \lim_{n\to\infty}d(y_n,\E_dT) = 0 \qquad \mbox{a.s.}\,.
\end{equation}
In other words, the set of empiric mean trees coincides with the singleton of the
$d$-mean if $n$ is large enough.
\end{thm}

When $\nu$ is an arbitrary probability measure on $\TT$, it may give positive
mass to sets of trees with infinitely many nodes. First we state the strong law
of large numbers for random elements taking values in a compact metric space
given in Sverdrup-Thygeson (1981). This covers the space of trees with infinite
number of vertices. Then we show that the metric space $\TT$ is compact; this
implies in particular that the expected tree is well defined ($\E_d(T)$ is non
empty).

Consider a compact metric space $(\KK,d)$. Let $\mc B$ denote the $\sigma$-field
generated by the open sets, and so the elements of $\mc B$ are the Borel sets.
Let $\nu$ be a probability measure on $\mc B$.  We define the expected value
with respect to the measure $\nu$ and the distance $d$ following the ideas
developed in the previous section.  Let $g\colon \KK\to \bb R^+$ be given by
\begin{equation}
  \label{px7}
  g(y):=\int_{\KK}d(y,x)\;\nu(dx).
\end{equation}
  Since
\begin{equation}
\label{px8}
\vert g(y)-g(t)\vert \leq \int_{\KK}\vert
d(y,x)-d(t,x)\vert\;\nu(dx) \leq\int_{\KK}d(y,t)\nu(dx)
=d(y,t)\;,
\end{equation}
we get that $g$ is Lipschitz continuous. Since it is defined on a compact space,
it attains its minimum. This shows that the $d$-mean set $\E_d(\nu)$ defined as in
\reff{p4} is non empty. The empiric mean $\overline T_n$ is defined as in
\reff{p5}.

\begin{thm}[Sverdrup-Thygeson, 1981]
Let  $\nu$ be a probability on the compact metric space  $(\KK,d)$ such that
$\E_d(\nu)$ has only one point. Consider
$\{T_n:n\geq 1\}$, an i.i.d. sample for $\nu$. Then, the empirical
$d$-centers converge uniformly to $\E_d(\nu)$ almost surely:
\begin{equation}
\lim_{n\to\infty}\sup_{a\in \overline
T_n}d\Big(a,\E_d(\nu)\Big)=0\quad \hbox{a.s.}\,.
\end{equation}
\end{thm}

The results of this section can be extended to the following family of
functions $g_p$ defined for $p\geq 1$ by
\[
g_p(y)\;=\;\int_{\KK}d(y,x)^p\;\nu(dx)\;.
\]

\section{Examples}
\label{s5}

\paragraph{\it Mode parameter}
Consider a finite space $\KK$ with the discrete distance given
by
\begin{equation}
\label{ddiscreta} d(x,y)\;=\; \left\{
\begin{array}{ll}
0&\mbox{if} \hspace{2mm} x=y \;,\\
1&\mbox{otherwise.}
\end{array}
\right .
\end{equation}
In this case,
\begin{equation}
\label{gdiscreta}
g(x)=\int_{\KK}d(x,y)\;\nu(dy)=\sum_{y\not=x}\nu(y)=1-\nu(x)\;.
\end{equation}
So, the $d$-center parameter   for $(\KK,d,\nu)$ is just the mode of $\nu$.

\paragraph{\it Mean and median parameters}
Consider $\KK = [0,1]^n \subset \R^n$, and $d(x,y) = ||x-y||^p$.  Let $\nu$ be
any probability measure on $\KK$. Then, if $p=2$ we have that the $d$-center
parameter is the usual expected value. For $n=1$ and $p=1$ we get the median,
and for $n>1$ the spatial median or multivariate $L_1-median$, see for instance
Haldane (1948) and Milasevic and Ducharme (1987).

\paragraph{\it Product Space}
We say that $(\KK, d,\nu)$ is a centered space if it has a unique
$d$-center.
We now prove that the product of a finite number of centered spaces is a
centered space.

\begin{lemma}
Let $(\KK_i,d_i,\nu_i)$ be   spaces  with unique $d$-centers
$C_i=E_{d_i}(\nu_i)$, for  $i=1,2$.  Then, if we consider the
product space $\KK=\KK_1\times\KK_2$ with
\begin{equation}
\label{dproducto} d(\hat x,\hat y)=d_1(x_1,y_1)+d_2(x_2,y_2)\,,
\end{equation}
for $\hat x=(x_1,x_2)\in \KK$ and the product measures
$\nu=\nu_1\times\nu_2$, we get that $(\KK,d,\nu)$ has also a unique
$d$-center
  $(C_1, C_2)$.
\end{lemma}

\begin{proof}
We need to prove that $(C_1,C_2)$ is the unique point minimizing
$g:\KK\to\R$. We get
\begin{equation}
\label{gproducto} g(\hat x)=\int _{\KK_1}\int
_{\KK_2}\Big(d_1(x_1,y_1)+d_2(x_2,y_2)\Big)
\nu_2(dy_2)\nu_1(dy_1)=g_1(x_1)+g_2(x_2)\;,
\end{equation}
where
\begin{equation}
\label{cadagi} g_i( x)=\int _{\KK_i}d_i(x,y)\nu_i(dy)\,.
\end{equation}
from where the result follows.
\end{proof}

\section{Invariance Principle}
\label{s6}

In this section we consider a sequence of independent identically
distributed random trees $(T_1,\dots,T_n)$ with empiric mean
$g_n(t)$ given by \reff{px1} and prove an invariance principle for
the centered process \[(\sqrt n (g_n(t)-g(t)), t\in\TT)\,,\] as
$n\to\infty$. The main tool is the following general result.

\begin{thm}[Ledoux and Talagrand (1991) pag 395--396]
  \label{lt}
  Let $\TT$ be a compact metric space and $C(\TT)$ be the separable Banach
  space of continuous functions on $\TT$ with the sup norm. Let
  $(\Omega,\FF,\P)$ be a probability space and $X:\Omega\to C(\TT)$ be a
  random element of $C(\TT)$ with $\E X(t) =0$ and $\E X(t)^2 <\infty$ for
all
  $t$ in $\TT$. Assume that $X$ is Lipschitz, that is, there exists a
positive
  random variable $M$ with $\E M^2<\infty$ such that
\begin{equation}
  \label{x1}
  |X(\omega,s)-X(\omega,t)|\,\le\, M(\omega)\,d(s,t)\,,
\end{equation}
for all $\omega\in\Omega$, $s,t\in\TT$. Assume there exists a probability
measure $\mu$ on $(\TT,d)$ such that
\begin{equation}
  \label{x2}
  \lim_{\delta\to0}\, \sup_{t \in \cal T}\,\int_0^\delta \Bigl[-\log[\mu(B(t,u))]\Bigr]^{1/2}\,du
  \;=\;0\,,
\end{equation}
where $B(t,u)$ is the ball centered at $t$ with radius $u$. (This is called
the majorizing measure condition for $(\TT,d)$.) Then $X$ verifies the Central
Limit Theorem in $C(\TT)$. That is, if $X_1,\dots,X_n$ are i.i.d. with the
same law as $X$, then $n^{-1/2}(X_1+\dots+X_n)$ converges to a Gaussian
process with mean zero and the same covariance function as $X$.
\end{thm}

\paragraph{\bf The majorizing measure condition} If $\TT$ is finite, the
condition is satisfied automatically by any measure $\mu$ on $\TT$ giving
positive mass to all elements of $\TT$. Indeed, $\mu(B(t,u))\ge\mu(t)>0$ and
the integral in \reff{x2} is dominated by $[-\log(\mu(t))]^{1/2}\delta$.

\begin{lemma}
  \label{y1}
  Let $\TT$ be the set of trees. Let $0<z<m^{-3/2}$ and $\phi$ defined by
\begin{equation}
  \label{x33}
  \phi(v) = z^{\gen(v)}\,.
\end{equation}
Then the majorizing measure condition is satisfied for $(\TT,d)$ with the
distance defined by \reff{p2} and  this $\phi$.
\end{lemma}

\proof Since for finite trees the result follows, we assume the trees in $\TT$
have infinitely many generations. Define the cylinder of generation $k$ induced
by the tree $t\in\TT$ by
\begin{equation}
  \label{x9}
  \TT_k(t):= \{s\in\TT\,:\, s(v)=t(v)\hbox{ if  }\gen(v)\le k\}\,.
\end{equation}
%
Define for $u\geq 0$
\begin{eqnarray}
  \label{x8}
  k(u)\;=\;k(u,\phi)
  &:=&\inf\Bigl\{k\,:\, \sum_v\phi(v)\one\{\gen(v)>k\}<u\Bigr\}\,.\nonumber
\end{eqnarray}
Since $\sum_v \phi(v)<\infty$, $k(u)$ goes to $\infty$ as $u$
goes to 0.
 Since
\begin{equation}
  \label{34}
  \sum_v\phi(v)\one\{\gen(v)>k\} = \sum_{i> k} m^{i-1} z^i =
\frac{z(mz)^k}{1-mz}\;,
\end{equation}
we can write
\begin{equation}
  \label{34a}
  k(u) =\inf\{k\,:\, z(mz)^k/(1-mz) <u\}\,.
\end{equation}
We have
\begin{equation}
  \label{x7}
\TT_{k(u)}(t)\subset  B(t,u)\,.
\end{equation}

A natural choice for a majorizing measure in $\TT$ is the measure
induced by the product measure $\nu_\rho$ on $\{0,1\}^{\widetilde V}$
with marginals $\nu_\rho \{\xi\,:\, \xi(v)=1\}=\rho$, for
$v\in\widetilde V$. Given a configuration $\xi\in\{0,1\}^{\widetilde
  V}$, define $x(\xi)$ as the maximal tree from the root whose
vertices are contained in the set $\xi$.  In other words, inductively, $x(\xi)(1)=\xi(1)$ and

\begin{equation}
x(\xi)(va)\;:=\;
\left\{
\begin{array}{ll}
1& \mbox{ if} \hspace{2mm} x(\xi)(v)=1 \hspace{2mm} \mbox{and} \hspace{2mm} \xi(va)=1\\
0& \mbox{ otherwise}\;,
\end{array}
\right.
\end{equation}
for each $v\in \{0,1\}^{\widetilde V}$ and $a\in A$.
Define the measure  $\mu_\rho$  induced on $\TT$ by this application:
$$
\mu_\rho(B):=\nu_\rho\{\xi:x(\xi)\in B\}\;.
$$
To check that $\mu_\rho$ is a majorizing measure, let $\beta>0$ be
defined by $e^{-\beta} =
\min\{\rho,1-\rho\}$. The number of vertices in the first $k$ generations of
the full tree is $(m^{k}-1)/(m-1)\le 2m^k$. Hence the probability of any
cylinder with $k$ generations is bigger than $e^{-2\beta m^k}$:
\begin{equation}
  \label{x31}
  \mu_\rho(\TT_k(t)) \ge\nu_\rho\Big\{\xi\in \{0,1\}^{\widetilde
V}:\xi(v)=t(v)\;\hbox{if $\gen(v)\leq k$}\Big\}
  \ge e^{-2\beta m^k}\,.
\end{equation}
uniformly in $t$. This and \reff{x7} imply that the supremum of the integral in
\reff{x2} is bounded above by
\begin{equation}
  \label{x32}
  \int_0^\delta (2\beta)^{1/2} m^{k(u)/2} du
\;=\;\int_0^\delta (2\beta)^{1/2} e^{k(u)\log (m^{1/2})} du
\;\le\; (2\beta)^{1/2} \int_0^\delta \frac{1}{u^{1-\eps}} du \,,
\end{equation}
for $\delta$ small enough, if there exists an $\eps>0$ such that $k(u)\le -(\log
u)(1- \eps)/\log(m^{1/2})$, for $u$ small enough.  In this case the proof is
finished because  for $\eps>0$ \reff{x32} converges to zero as $\delta\to 0$.
Call $\gamma=(1- \epsilon)/\log(m^{1/2})$. In view of \reff{34},
we look for $\gamma>0$ such that $ (mz)^{-\gamma\log u} <
u(1-mz)/z $. That is,
\[
u^{-\gamma\log (mz)-1} < \frac{1-mz}{z}\;.
\]
For $u$ sufficiently small it suffices that $-\gamma\log (mz)-1 >0$ and
$z<m^{-1}$.  Substituting $\gamma$ and noticing that $\log (mz)<0$, we need to
find an $\eps>0$ such that
\[
-(1-\eps) < \frac{\log(m^{1/2})}{\log (mz)}, \quad\hbox{ that is,
}\quad \eps <1+ \frac{\log(m^{1/2})}{\log (mz)}\;,
\]
which exists since $z <m^{-3/2}$. \endproof

We are now able to obtain the asymptotic distribution of the process
$$
\sqrt n \left(g_n(t)-g(t)\right) = \frac{
\sum_{i=1}^n[d(T_i,t)-\E(d(T_i,t))]}{\sqrt n}.
$$

\begin{thm}
\label{test}
Let $\TT$ be the set of trees with at most $m$ offspring. Consider the distance
given in \reff{p2} for $\phi(v) =z^{\gen(v)}$ with $0<z<m^{-3/2}$. Let $\{ T_i :
i \ge 1\}$ be a sequence of i.i.d.  random trees on $\TT$ with the same law as
$T$. Then the process $(\sqrt n (g_n(t)-g(t)), t\in\TT)$ converges weakly as
$n\to\infty$ to a Gaussian process $W$ with zero mean and the same covariance
function as the process $X\in(\R^+)^\TT$ defined by $X(t) = d(T,t) -
\E(d(T,t))$.
\end{thm}

\begin{proof}
Since $\big|\, X(\omega,t) - X(\omega,t') \big| \le 2 d(t,t')$
the result follows from the previous Lemma and Theorem (\ref{lt}).
\end{proof}

\section{Statistical applications}
\label{s7}

Let $T$ be a random tree in $\TT$ with distribution $\nu$ and mean
distances $(g(y),\,y\in\TT)$ defined in \reff{px10}. Let $\nu_0$ be a
distribution on the tree space $\TT$ with mean distances
$(g_0(y),\,y\in\TT)$. The goal is to test

\noindent H0: $\nu = \nu_0$\hfill\break
HA: $\nu \neq \nu_0$

\noindent using an i.i.d. sample of random trees $\{T_i: i\ge
1\}$. Notice however that the rejection of H0 does not imply the
rejection of $\E T = \E_d(\nu_0)$.

To perform the test we propose the statistic
\begin{equation}
  \label{p29}
  \sup_{y\in \TT} |W_n(y)| = \sup_{y\in \TT} \sqrt n\,
\big|\,g_n(y)-g_0(y))\,\big|,
\end{equation}
whose asymptotic law under H0 is obtained from Theorem
(\ref{test}) and the Continuous Mapping Theorem.  We reject the
null hypothesis at level $\alpha$ if
$$
\sup_{y\in \TT} |W_n(y)|>q_\alpha\;,
$$
where $q_\alpha$ satisfies $P(\sup_{y\in \TT}
|W(y)|>q_\alpha)=\alpha$, for $W$ given in Theorem \reff{test} under
$\nu=\nu_0$.

The test rejects $\nu=\nu_0$ if $g$ determines $\nu$ unequivocally.

In practice the distribution of $\sup_{y\in \TT} |W(y)|$ depends on
the covariance of the process $X(t) = d(T,t) - \E(d(T,t))$ which in
general is unknown. A possible way to deal with this problem is to
approximate $q_\alpha$ using bootstrap.  The validity of the bootstrap
in this context remains an open problem. Alternatively, one can
simulate trees with distribution $\nu_0$ and estimate $q_\alpha$.

For the problem of two samples (of same size, for instance) one may
use the statistic
\begin{equation}
  \label{p30}
  \sup_{y\in \TT} \sqrt n|g_n(y)-g'_n(y)|\,,
\end{equation}
where $g_n$ and $g'_n$ correspond to the samples of $T$ and $T'$
respectively.

\paragraph{When  $g$ characterizes the measure $\nu$?}

Busch et al (2006) prove that $g=(g(t),t\in\TT)$ characterizes the
vertex-marginal distributions as follows. Let $\nu$ and $\nu'$ be two measures
in $\TT$ and $g$, $g'$ be the corresponding processes. Then $g=g'$ if and only
if $\nu\{t:t(v)=1\} = \nu'\{t:t(v)=1\}$ for all vertex $v$. In that paper it is
proven that under certain Markov hypothesis, the vertex-marginals identify
univoquely the measure. The class of random trees satisfying those hypothesis
includes Galton-Watson branching processes and other related processes.

\section{Metrics and negative curvature}
\label{s8} In this section we show that our tree space cannot be
embedded in a metric space of non positive curvature. Then we
discuss other possible metrics that have been considered for
spaces of trees.  A natural way of embedding the discrete tree
space $\TT$ in a continuous space would be to consider a tree as a
function $x:\widetilde V\to\R^+$ (instead of $\{0,1\}$), where the
value $x(v)$ would represent the length of the edge connecting the
node $v$ to her mother. The value $x(v)=0$ means that the node $v$
is not present. The metric could be the one given in \reff{p2}
which coincides with the previous one for trees with unitary edge
lengths. A tree condition like ``$x(va)>0$ implies $x(v)>0$'' is
also needed, but other conditions could be proposed. For instance
one could collapse the vertices with $x(v)=0$ but in this case the
trees would not have a limited number of offspring nor the vertex
notation introduced in Section \ref{s1} would be appropriate.

Let $a,b,c,x$ be arbitrary distinct points in a metric space $\TT$
such that $x$ belongs to a geodesic from $a$ to $b$, that is,
$d(a,b)=d(a,x)+d(x,b)$. Let $a',b',c',x'$ in $\R^2$ be points located
in such a way that the relative distances are the same, that is,
$d(a,b)=d'(a',b')$, $d(a,c)=d'(a',c')$, etc, where $d'$ is the
Euclidean distance in $\R^2$. It is said that $\TT$ is of \emph{non
  positive curvature} if $d(x,c)\le d'(x',c')$ for any choice of
$a,b,c,x$. These spaces are also called {\sl CAT}$(0)$, see Billera, Holmes
and Vogtmann (2001), page 750.

We now give an example showing that our space cannot fit the above property.
Let $a=\{111,12\}$ and $b=\{11,121\}$ be the trees in Figure 1 and $x=\{11,12\}$ and
$c=\{111,112,121,122\}$ those of Figure 2.
Consider the distance \reff{p2} with $\phi(v)$ depending only on the
generation of $v$, so that $\phi(111)=\phi(121)=\phi(122)=\phi(112)=\alpha$,
for some $\alpha>0$.  The tree $x$ belongs to a geodesic between $a$ and $b$:
$d(a,x) = d(b,x)=\alpha$, $d(a,b)= 2\alpha$. On the other hand $d(a,c) =
d(b,c)= 3\alpha$ and $d(x,c) =4\alpha$. Consider the corresponding Euclidean
triangle $(a',b',c')$ with the same relative distances. The point equidistant
from $a'$ and $b'$ in the Euclidian geodesic, corresponding to $x$, is $x'=
(a'+b')/2$. Since $d'(x',c')= \sqrt 8 \alpha <4 \alpha=d(x,c)$, our tree space
cannot be embedded in a {\sl CAT}$(0)$ space.


\begin{figure}[ht]
\setlength{\unitlength}{0.5mm}
\phantom{1}\hskip 2cm
\begin{picture}(200,100)
\thicklines
\put(100,90){$(c)$}
\put(120,50){\circle*{5}}\put(118,35){1}
\put(120,50){\vector(4,3){38}}
\put(120,50){\vector(4,-3){38}}
\put(160,80){\circle*{5}}\put(158,65){12}
\put(160,20){\circle*{5}}\put(158,5){11}
\put(160,80){\vector(2,1){38}}
\put(160,80){\vector(2,-1){38}}
\put(200,100){\circle*{5}}\put(198,90){122}
\put(200,60){\circle*{5}}\put(198,65){121}
\put(200,40){\circle*{5}}\put(198,28){112}
\put(200,0){\circle*{5}}\put(198,5){111}
\put(160,20){\vector(2,1){38}}
\put(160,20){\vector(2,-1){38}}
\put(-20,90){$(x)$}
\put(0,50){\circle*{5}}\put(-2,35){1}
\put(0,50){\vector(4,3){38}}
\put(0,50){\vector(4,-3){38}}
\put(40,80){\circle*{5}}\put(38,65){12}
\put(40,20){\circle*{5}}\put(38,5){11}
\end{picture}
\caption{The trees $x=\{11,12\}$ and
c=\{111,112,121,122\}.}\label{fig2}
\end{figure}
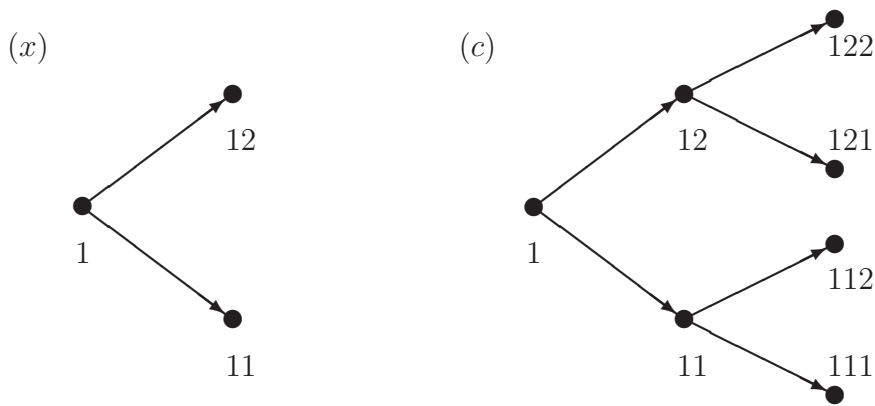

Notice that the tree $\bar x:=(111,121)$ is also in a (different)
geodesic between $a$ and $b$. However $d(\bar x,c) = 2\alpha<\sqrt
8\alpha=d(x,c)$. In fact the triangle $(a',b',c')$ is the same and $\bar x'=
x'$. Actually, the fact that there are two geodesics going from $a$ to $b$
indicates that the space cannot be of negative curvature.

The metric we propose for the space of trees is usual in interacting particle
systems, some of which are defined in $\{0,1\}^S$ for countable $S$, for
example. The product of the discrete topologies induces a metric like
\reff{p2}. Under this metric, the convergence of a sequence $x_n$ to $x$ is
equivalent to the convergence of $x_n(v)$ to $x(v)$ for all vertex $v$.
Valiente (2001) considers spaces of finite trees with ordered
vertices, reviews several distances and proposes a new metric. An important
example is the so called ``edit distance'', which counts the number of
operations (eliminate a vertex, add a vertex) that need to be done in order to
transform one tree into another one. Critchlow (1980) proposes some metrics in
the set of permutations of a finite sequence that may be adapted to a finite
space of trees. It would be nice to understand if our results can be proven in
those spaces.

Billera, Holmes and Vogtmann (2001) describe various spaces of
``phylogenetic trees'' and construct a (continuous) convex metric
space of trees with a fixed number $n$ of final vertices (i.e.,
vertices with no daughters). The resulting space ${\cal F}_n$ is
{\sl CAT}$(0)$. Phylogenetic trees are constructed from the final
vertices to the root by successively grouping subsets of vertices
as in the Kingman's coalescent. Each vertex with descendants
represents the most recent common ancestor of the descendants, and
the length of the edge $(v,v')$ represent the time a group of
species represented by $v'$ needed to split.  In our space the
trees can have variable number of final vertices; our
counterexample does not apply to spaces of trees with fixed number
of final vertices. Another difference with phylogenetic trees is
that in our space we do not label the (final) vertices.

\section{Final remarks}

Our motivation was to produce a statistical tool to study the asymptotic
behavior of sequences of random trees. The law of large numbers is not directly
applied to construct the tests, but is important to guarantee the consistency of
the estimators.   On the other hand, the central limit theorem (Theorem
\ref{test}) uses the tree structure and a particular form of the distance. The
shape of the function $\phi$ intervening in the distance was necessary to show
that the majorizing measure condition holds (Lemma \ref{y1}). We believe this
can be extended to other structures contained in a subset of $\{0,1\}^S$ for $S$
countable. Another possible extension is to eliminate the upperbound $m$ on the
number of offsprings. If the mean number of offsprings is not finite, then the
limits may be stable laws, but this is to be stablished.

The statistical application we have considered in Section 7 points
in the direction of a Kolmogorov-Smirnov type goodness of fit
test. We are interested in the decision problem: given a random
sample $T_1,...,T_n$ can we decide if their underlying common
distribution $P$ is a given $P_0$? For instance does the sample
follow the Galton-Watson model with parameter $p_0$? We think that
the statistic given in section 7 is adequate for this problem. The
results in Busch et al (2006) where our test has been applied to
several simulated examples, and a real data example to classify
FGF protein families points in this direction.

The implementation of the tests requires the computation of the
statistic \reff{p29} which is a supremum over the space of trees
of the distance of the tree to the mean tree.  The computation
time of this task may increase fast with the number of nodes.
Busch et al (2006) propose a method to transform this problem in
the computation of the minimal cut of the flux of a related graph.
This allows to see the behavior of the test in some concrete
examples related to Galton-Watson generated random trees.

\section*{Acknowledgements}
We thank one of the referees for mentioning that our setup was
proposed by Otter and Neveu. We thank all referee´s for their
careful reading and constructive comments.

This work started in the open problem session of the 6th Brazilian School of
Probability, Ubatuba, 2002. This paper is partially supported by the Institute
of Millenium for the Global Advance of Mathematics, Funda\c c\~ao de Amparo
\`a Pesquisa do Estado de S\~ao Paulo, Conselho Nacional de Desenvolvimento
Cient\'\i fico e Tecnol\'ogico, Programa N\'ucleos de Excel\^encia and the
Brazil-Argentina Agreement CAPES-SECyT.


\bigskip
\begin{thebibliography}{60}

\bibitem{[BHV]} Billera, L. J.; Holmes, S. P.; Vogtmann, K. (2001) Geometry of
  the Space of Phylogenetic Trees, {\sl Advances in Applied Mathematics, \bf
    27}, 4:733-767.


\bibitem{bfffgl}   Busch, J. R.;  Ferrari, P. A.;  Flesia, A. G.; Fraiman, R.;
  Grynberg,  S. ; Leonardi, F. (2006)
Testing statistical hypothesis on random trees. Preprint.


\bibitem{[dF]} De Fitte, P. R. (1997) Th\'eor\`eme ergodique ponctuel
  et lois fortes des grandes nombres pour des points al\'eatoires d'un
  espace m\'etrique \`a courbure n\'egative. {\sl Ann.\ Probab. \bf 25}
  2:738-766.

\bibitem{[eh]} Es-Sahib, A.; Heinich, H. (1999) Barycentre canonique
  pour un espace m\'etrique \`a courbure n\'egative. {\sl S\'eminaire de
  Probabilit\'es, XXXIII, 355--370, Lecture Notes in Math., \bf
    1709}, Springer, Berlin.


\bibitem{[Haldane]} Haldane, J. B. S. (1948). Note on the median of a
  multivariate distribution.  {\sl Biometrika \bf 35}, 414-415.

\bibitem{[H]} Herer, W. (1992) Mathematical expectation and strong
  law of large numbers for random variables with values in a metric
  space of negative curvature. {\sl Probab.\ Math.\ Statist. \bf 13}
  1:59-70.

\bibitem{[K]} Kingman, J.F.C. (1982) The coalescent. {\sl Stoch.\
    Proc.\ Appl., \bf 13} 235-248.


\bibitem{[KM]} Kurata, K. and Minami, N. (2004). The equivalence of two constructions of
Galton--Watson processes. Preprint mp\_arc 04-295.


\bibitem{[lt]} Ledoux, M. and Talagrand, M. (1991) {\it Probability in
    Banach Spaces. Isoperimetry and Processes.}  Springer-Verlag.


\bibitem{Liggettbook} Liggett, T.~M., (1985) \emph{Interacting Particle
    Systems}.  Springer-Verlag, New York.

\bibitem{[MD]} Milasevic, P. and Ducharme, G. R. (1987). Uniqueness of the spatial median.
{\sl Annals of Statistics \bf 15}, 1332-1333.

\bibitem{[M]} Mourier, E. (1953) El\'ements al\'eatoires dans un
  espace de Banach. {\sl Ann.\ Inst.\ Henri Poincar\'e, \bf 13}
  159-244.

\bibitem{Neveu}
Neveu, J. (1986). Arbres et processus de Galton-Watson.
{\sl Ann. Inst. H. Poincar\'e, Probabilit\'es et
Statistique} {\bf 22,} 199--207.

\bibitem{Otter}
Otter, R. (1949).  The Multiplicative Process.
{\sl Ann. Math. Stat.} {\bf 20,} 206--224.

\bibitem{[S]} Sverdrup-Thygeson, H. (1981). Strong law of large numbers for
measures of central tendency and dispersion of random variables in compact
metric spaces. {\sl Ann. \ Statist., \bf 9} 141-145.

\bibitem{[V]} Valiente, G. (2001) An Efficient Bottom-Up Distance
  between Trees.  {\sl Proc.\ 8th Int.\ Symposium on String Processing
    and Information Retrieval}, IEEE Computer Science Press. SPIRE
  2001: 212-219.
\end{thebibliography}
\end{document}